\documentclass[12pt,a4paper]{article}
      \usepackage[english]{babel}
      \usepackage{amsmath}
      \usepackage{amstext}
      \usepackage{amssymb}
      \usepackage{amsthm}
      \usepackage{amsfonts}
      \usepackage{fullpage}
      \usepackage{graphicx}
      \usepackage{epsf}
      \usepackage[dvips]{epsfig}
      \usepackage{stmaryrd}
      \usepackage[T1]{fontenc}
      \usepackage{latexsym}
      \usepackage{psfrag}
      \newtheorem{prop}{Proposition}[]
      
      \newenvironment{demo}{\noindent {\bf Proof:}}{\hfill $\square$ \\}
      \newtheorem{thm}{Theorem}[]
      \newtheorem{cor}{Corollary}[thm]

      \newcommand{\R}{\mathbb R}
      \newcommand{\C}{\mathbb C}
      \newcommand{\Z}{\mathbb Z}
      \newcommand{\Q}{\mathbb Q}
      \newcommand{\defi}{\stackrel{\mbox {\tiny def}}{=}}
      
      \newenvironment{appli}{\left( \begin{array}{ccc}}{\end{array} \right)}
      \newcommand{\dans}{& \longrightarrow &}
      \newcommand{\donne}{& \longmapsto &}
      \newcommand{\ba}{\begin{appli}}
      \newcommand{\ea}{\end{appli}}

      \newcommand{\tq}{\ | \ }
      \newcommand{\B}{\mathbb B} 
      \renewcommand{\H}{\mathcal H}

      \renewcommand{\P}{{\mathcal{P}}}
      \renewcommand{\S}{{\mathcal{S}}}
      
      \renewcommand{\C}{{\mathcal{C}}}
      \renewcommand{\B}{{\mathcal{B}}}

\title{A homological condition for a dynamical and illuminatory classification of torus branched coverings}
\author{Thierry Monteil \footnote{Institut de Math\'ematiques de Luminy, CNRS UMR 6206,
Case 907, 163 Avenue de Luminy, 13288 Marseille cedex 09, France
-- E-Mail: monteil@iml.univ-mrs.fr -- Tel: +33 4 91 26 96 77 --
Fax : +33 4 91 26 96 55 --}}
\date{}

\begin{document}

\maketitle


\vspace{2 cm}

\begin{abstract} 

\noindent We prove that, for translation surfaces whose homology is generated by the periodic orbits,
the notions of 
\begin{itemize}
  \item finite blocking property
  \item pure periodicity
  \item torus branched covering
\end{itemize}
are equivalent. 
In particular, we prove this equivalence for convex surfaces 
and on a dense open subset of full measure on every normalized stratum.\\


{\em \noindent AMS classification: 
37D50, 
37C27, 
57M12, 
51E21, 
55N99, 
37E35, 
32G15. 

}

\end{abstract}

\clearpage

\section*{Introduction}

A {\em translation surface} is a triple $(\S, \Sigma, \omega)$
such that $\S$ is a topological compact connected surface,
$\Sigma$ is a finite subset of $\S$ (whose elements are called
{\em singularities}) and $\omega = (U_i,\phi_i)_{i\in I}$ is an
atlas of $\S \setminus \Sigma$ (consistent with the topological
structure on $\S$) such that the transition maps (i.e. the $\phi_j
\circ \phi_i^{-1} : \phi_i(U_i\cap U_j) \rightarrow \phi_j(U_i\cap
U_j)$ for $(i,j)\in I^2$) are translations. This atlas gives to
$\S \setminus \Sigma$ a Riemannian structure; therefore, we have
notions of length, angle, measure, geodesic... Moreover, we assume
that $\S$ is the completion of $\S \setminus \Sigma$ for this
metric.

A translation surface can also be seen as a holomorphic
differential $h$ on a Riemann surface (the singularities
correspond to the zeroes of the differential, and in an admissible
atlas $\omega$, $h$ is of the form $h=dz$). 

Translation surfaces provide one of the main tool for the study of rational polygonal billiards.\\

Since the unit tangent bundle of $\S$ enjoys a canonical global
decomposition $U\S=\S \times \mathbb{S}^1$ (the rotational holonomy is trivial), 
the study of the geodesic flow on $\S$ can be done
through two points of view depending on whether the variable is
the first or the second projection:

\begin{description}

  \item[Dynamics] We can fix one particular direction $\theta\in\mathbb{S}^1$.

    This corresponds to the study of the {\em directional flow} $\phi_{\theta} : \S\times \R \rightarrow \S $.

    In that context, we say that a translation surface $\S$ is {\em purely periodic} if
    for any $\theta\in\mathbb{S}^1$, the existence of a periodic orbit
    in the direction $\theta$ implies that the directional flow $\phi_{\theta}$ is periodic 
    (i.e. there exists $T>0$ such that $\phi_{\theta}^T=Id_{\S}$ a.e.).

  \item[Illumination] We can also fix one point in $x\in\S$. 

    This corresponds to the study of the {\em exponential flow} $exp_x : \mathbb{S}^1 \times \R \rightarrow \S $.

    In that context, we say that a translation surface $\S$ has the {\em finite blocking property} 
    if for every pair  $(O,A)$ of points in $\S$,
    there exists a finite number of points $B_1, \dots , B_n$
    (different from $O$ and $A$) such that every geodesic from $O$ to $A$ meets one of the $B_i$'s.

\end{description}

A {\em torus branched covering} is a translation surface $\S$ such that there exists a branched translation covering
from $\S$ to a flat torus (a {\em branched translation covering} between two translation surfaces is a
map $\pi : (\S,\Sigma) \rightarrow (\S',\Sigma')$ that is a topological branched covering that locally preserves the
translation structure.).
This is a global geometric property.\\

In \cite{Monteil-Veech} and \cite{Monteil-pp}, the three notions of torus branched covering, finite blocking property and pure periodicity
have been proved to be equivalent for Veech surfaces and surfaces of genus two.
Gutkin proved independently that the notions of torus branched covering and finite blocking property 
are equivalent for Veech surfaces \cite{Gutkin}.

In this paper, we prove that the equivalence is true for any surface whose homology is generated by the periodic orbits of 
the geodesic flow.

In particular, the three notions are equivalent for convex surfaces and on a dense open subset of full measure of any stratum.

The paper is organized as follows :
in section \ref{section-background}, we recall some background and useful tools.
The second one is devoted to the proof of the following results :\\

{\bf Theorem \ref{th-hol=2}.} {\em
Let $\S$ be a purely periodic translation surface. 

Then $hol(\P(\S))$ is a lattice of $\R^2$.}\\

{\bf Theorem \ref{th-equivalence-homologie}.} {\em
Let $\S$ be a translation surface whose homology is generated by the
periodic orbits of the geodesic flow on $\S$.
Then the following assertions are equivalent: 
\begin{enumerate}
  \item $\S$ is a torus branched covering
  \item $\S$ has the finite blocking property
  \item $\S$ is purely periodic
\end{enumerate} }

In section \ref{section-generation-convex}, we try to see to what extent does this result applies.
In particular, we prove\\

{\bf Proposition \ref{th-convex}.} {\em
Let $\S$ be a convex surface, then the periodic orbits of $\S$ generate $H_1(\S,\Z)$.}\\

{\bf Proposition \ref{th-generique}.} {\em
In any stratum $\mathcal{H}_1(k_1,\dots,k_n)$, 
the set of translation surfaces whose homology is generated by the perodic orbits 
is an open dense subset of full measure.}\\

The last section is a discussion about the tools we used, 
particularly the $J$-invariant (in the proof of Theorem \ref{th-hol=2}) 
and the existence of two transverse parabolic elements 
(for the proof of the result about Veech surfaces), 
and the fact that in general one needs to find periodic orbits in more than two directions 
to generate the homology of a translation surface.
We take the opportunity to introduce the class of translation surface whose $J$-invariant is of the simplest form, 
and present a collaborative encyclopedy of concrete translation surfaces.

\clearpage

\section{Background and tools}\label{section-background}

\subsection{Cylinder decomposition}

Let $\S$ be a translation surface.
A {\em cylinder} $\C$ of $\S$ is a maximal isometric copy of $\R / w \Z \times ]0,h[$ in $\S$ ($w>0$, $h>0$). 
The parameters $w$ and $h$ are unique and called the {\em width} (or {\em perimeter}) and the {\em height} of $\C$.
The {\em direction} of $\C$ is the direction of the image of $\R / w \Z \times \{ h /2 \}$
(which is considered here as an oriented closed geodesic); it is defined modulo $2\pi$.\\

We say that $\S$ admits a {\em cylinder decomposition} in the direction $\theta$ 
if the (necessarily finite) union of cylinders in that direction is dense in $\S$.
The remaining is a finite union of saddle connections 
(a {\em saddle connections} is a geodesic $\gamma: [0,1] \rightarrow \S$ joining two singularities 
and such that $\gamma(]0,1[)$ does not contain any singularity).\\

Any periodic trajectory for $\phi_\theta$ can be thickened to obtain a cylinder in the direction $\theta$.
Thus, saying that $\S$ is purely periodic is equivalent to say that in each direction $\theta$ that admits a periodic orbit, 
$\S$ admits a cylinder decomposition whose cylinders have commensurable perimeters. 
This property is stronger than the notion of complete periodicity introduced by Calta \cite{Calta}.

\subsection{(Translational) holonomy}

Let $\S$ be a translation surface.
Let $\gamma:[0,1] \rightarrow \S$ be a continuous curve. 
Thanks to local charts, $\gamma$ can be lifted on $\R^2$ (or $\mathbb{C}$) 
to a curve $\overline{\gamma}$ defined up to translation.
Hence, $hol(\gamma) \defi \overline{\gamma}(1)- \overline{\gamma}(0)$ is well defined and is called the {\em holonomy} of $\gamma$.
If $h$ is the holomorphic form that defines $\S$,
we have $hol(\gamma) = \int_\gamma h$.\\

If $\gamma'$ is homologous to $\gamma$, then $hol(\gamma')=hol(\gamma)$.
If we are looking to closed curves, $hol$ can be extended by formal sum to a morphism $hol:H_1(\S,\Z) \rightarrow \R^2$.

\mathversion{bold}
\subsection{Moduli space and $SL(2,\R)$-action}
\mathversion{normal}

A singularity $\sigma \in \Sigma$ has a conical angle of the form
$2(k+1)\pi$, with $k\geq 0$; we say that $\sigma$ is of {\em
multiplicity} $k$. 
In terms of holomorphic forms, it is equivalent to say that
there is a chart around $\sigma$ such that $h=z^{k}dz$ (i.e. $\sigma$ is a zero of order $k$ of $h$).\\

If $1\leq k_1\leq k_2 \leq \dots \leq k_n$ is a sequence of positive integers whose sum is even, 
we denote by $\H(k_1, k_2, \dots, k_n)$ the {\em stratum} of translation surfaces with
exactly $n$ singularities whose multiplicities are $k_1, k_2, \dots, k_n$ 
(we consider only surfaces without {\em removable} singularities i.e. singularities of multiplicity $0$). 
A translation surface in $\H(k_1, k_2, \dots, k_n)$ has genus $g=1+(k_1+k_2+ \dots +k_n)/2$.\\

For any translation surface $\S$ and any $A\in SL(2,\R)$, we can define the translation surface 
$$A .(\S, \Sigma,(U_i,\phi_i)_{i\in I}) \defi (\S, \Sigma, (U_i,A \circ \phi_i)_{i\in I})$$ 
hence we have an action of $SL(2,\R)$ on the moduli space of translation surfaces. \\

Each stratum carries a natural topology and a $SL(2,\R)$-invariant
measure that are for example defined in \cite{Kontsevich}.
Let $\S$ be an element of some $\H(k_1, k_2, \dots, k_n)$ and 
let $\B$ be a basis of the {\em relative} homology of $\S$:
it is just the concatenation of a basis of the first topological homology group $H_1(\S,\Z)$ 
with a set of $n-1$ curves from a singularity to the other ones.
If $\S'$ is another translation surface (built on the same topological surface), 
let us denote by $hol_{\S'}$ the associated holonomy.
The map 
$$\Phi \defi \ba \H(k_1, k_2,\dots, k_n)  \dans \mathbb{C}^{2g+n-1} \\ 
\S'  \donne (hol_{\S'}(\gamma_1) , \dots , hol_{\S'}(\gamma_{2g+n-1}))  \ea$$
is named the {\em period map} and is a local homeomorphism in a neighbourhood of $\S$ 
and in this system of coordinates, the former measure is absolutely continuous relatively to Lebesgue.\\

In each stratum $\H(k_1, k_2, \dots, k_n)$, 
let  $\H_1(k_1, k_2, \dots, k_n)$ denotes the real hypersurface defined by the equation $area(\S)=1$.
The topology can be induced and the measure can be pushed from $\H(k_1, k_2, \dots, k_n)$ 
to $\H_1(k_1, k_2, \dots, k_n)$ that is stable under $SL(2,\R)$.

Masur \cite{Masur-ei} and Veech \cite{Veech-ei} proved that the volume of any such (normalized) stratum is finite, 
and that the action of $SL(2,\R)$ is ergodic on any connected component of any normalized stratum.

\mathversion{bold}
\subsection{The $J$-invariant}
\mathversion{normal}

In \cite{KenyonSmillie}, Kenyon and Smillie define an algebraic invariant for translation surfaces: the $J$-invariant. 
It takes values in the alternating product of $\R^2$ by itself
denoted here by $\R^2 \wedge_{\Q} \R^2$
(note that $\R^2$ is considered here as an infinite dimensional $\Q$-vector space).
If $(\beta,\leq)$ is a totally ordered basis of $\R^2$, 
then $\{u\wedge v \tq (u,v)\in \beta^2 \mbox{ and } u<v \}$ is a basis of $\R^2 \wedge_{\Q} \R^2$.\\

First, they define it for planar polygons: 
if $P$ is a polygon of $\R^2$ with vertices $v_1, \dots , v_n$ in counterclockwise order around the boundary of $P$, 
then $$J(P)\defi v_1 \wedge v_2 + v_2 \wedge v_3 + \dots + v_n \wedge v_1$$
If $P$ is a parallelogram, we have $J(P)=2 e_1 \wedge e_2$, 
where $e_1 = v_2-v_1$ and $e_2=v_3-v_2=v_4-v_1$ are two consecutive edges of $P$ (for the same counterclockwise order).\\

Then, if a translation surface $\S$ admits a cellular decomposition into planar polygons $\S=P_1 \cup \dots \cup P_n$, 
they define $$J(\S) \defi J(P_1)+\dots+J(P_n)$$
Any translation surface admits such a decomposition and 
$J(\S)$ is independent of the choice of the decomposition of $\S$ into
polygons (see \cite{KenyonSmillie}). 
Note that the singularities of $\S$ have to be located at the vertices of the $P_i$'s.\\

We have to notice that for any translation surface $\S$, $J(\S)$ cannot be equal to zero since 
if $\varphi$ denotes the linear map from $\R^2 \wedge_{\Q} \R^2$ to $\R$ such that for any $u$ and $v$ in $\R^2$
$\varphi(u \wedge v)=det(u,v)$, then $\varphi (J(\S)) = 2.area(\S) \neq 0$.

\section{A homological condition for a purely periodic translation surface to be a torus branched covering}

Thanks to Hopf-Rinow theorem, since $\S$ is a complete surface, 
any closed curve $\gamma$ on $\S$ can be tightened to give a closed geodesic in the same homology class. 
Therefore, $H_1(\S,\Z)$ is generated by the closed geodesics of $\S$.
Among them, are the periodic orbits for the geodesic flow. 
But they are not the only ones since there are also closed unions of ``consecutive'' saddle connections.
Let $P(\S)$ denotes the set of periodic orbits, 
and let $\P(\S)$ denotes the subgroup of $H_1(\S,\Z)$ generated by $P(\S)$.

\begin{thm}\label{th-hol=2}
Let $\S$ be a purely periodic translation surface. 

Then $hol(\P(\S))$ is a lattice of $\R^2$.
\end{thm}

\begin{demo}
\begin{description}
\item[Step 1: ] We prove that $dim_{\Q} (Vect_\Q(hol(P(\S))))=2$.

Let $V$ and $W$ be two non-colinear elements of $hol(P(\S))$. 
Up to let $SL(2,\R)$ act on $\S$, we can suppose that $V$ is horizontal whereas $W$ is vertical.

Now, assume by contradiction that there is an element $V'$ in $hol(P(\S))$ 
that is not in the $\Q$-vector sapce generated by $V$ and $W$.

$V'$ can be neither parallel to $V$ nor $W$ since $\S$ is purely periodic.

We write the real combination $V' = \alpha V + \beta W$.

We have supposed that $\alpha$ or $\beta$ is in $\R \setminus \Q$.
There is no loss of generality to assume that both $\alpha$ and $\beta$ are positive and that $\alpha$ is irrational
(see Figure \ref{3vecteurs.eps}, the positivity of $\alpha$ and $\beta$ is just here 
to avoid confusions in the computation of $J(\S)$).

\begin{figure}[h!]
  \begin{center}
    \psfrag{v}{$V$}
    \psfrag{v'}{$V'$}
    \psfrag{w}{$W$}
    \psfrag{a}{$\alpha$}
    \includegraphics[width=.43 \linewidth]{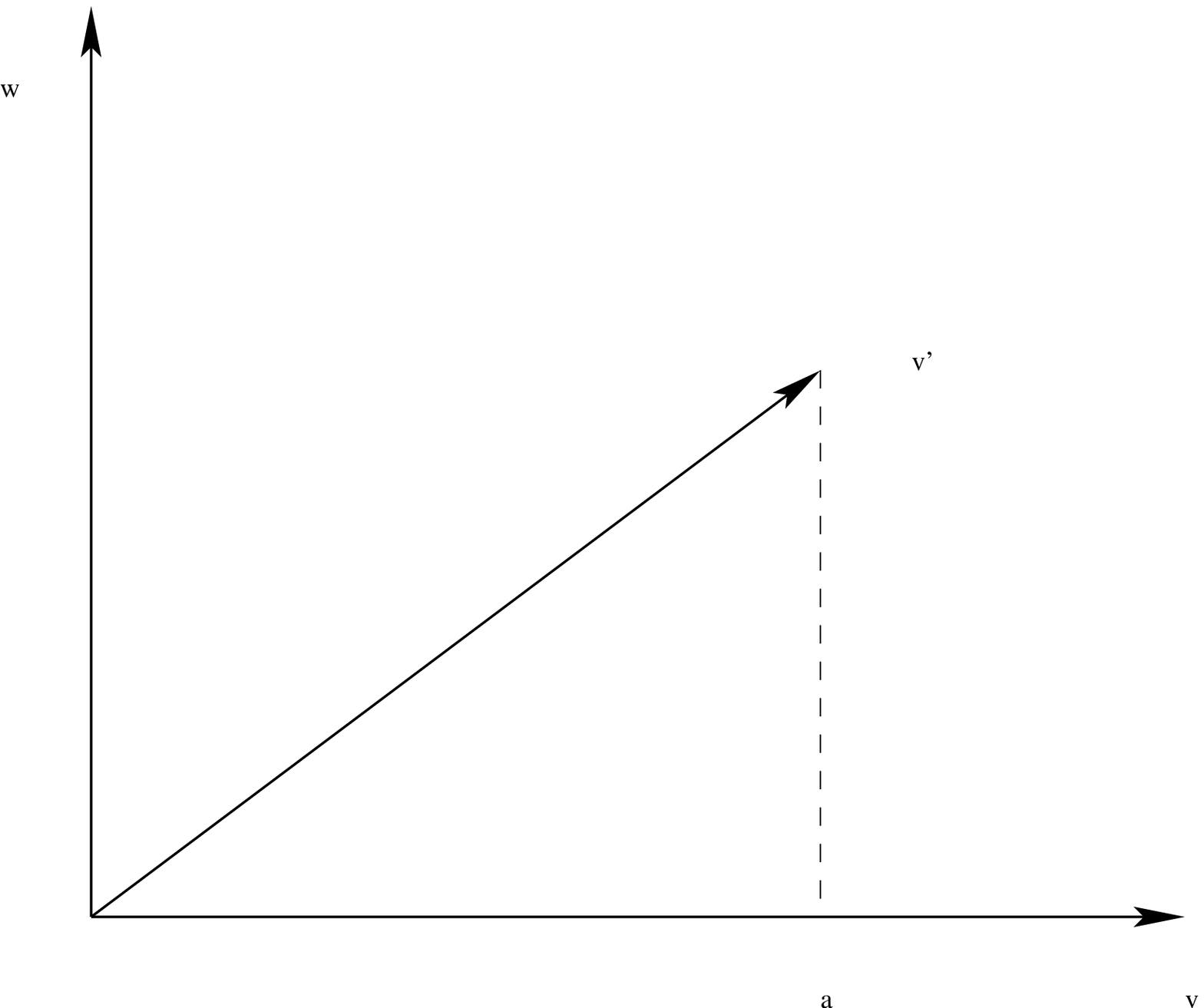}
    \caption{\label{3vecteurs.eps} Three vectors in $hol(P(\S))$ that are free over $\Q$.}
  \end{center}
\end{figure}

We will now calculate $J(\S)$ in two different manners, along the couple of directions $(V,W)$, and then along $(V',W)$.

Since $\S$ is purely periodic, we can decompose $\S$ into parallel horizontal and vertical cylinders 
$(\mathcal{C}_i)_{i=0}^p$ and $(\mathcal{D}_j)_{j=0}^q$.

For $i\leq p$ and $j \leq q$, $\mathcal{C}_i \cap  \mathcal{D}_j$ 
is the reunion of a familly of disjoint open parallelograms (rectangles) $(P_{i,j,k})_{k=0}^{r(i,j)}$.
Let us denote $v_{i,j,k}$ and $w_{i,j,k}$ for the sides of $P_{i,j,k}$ along directions $V$ and $W$.

Note that $w_{i,j,k}$ does not depend on $j$ nor $k$ since it represents the height of $\mathcal{C}_i$ along the direction $W$, 
so we can write $w_{i,j,k}=w_i$. Hence, $J(P_{i,j,k}) = 2 v_{i,j,k} \wedge w_i$.

Since $\S$ is purely periodic, for any $i\leq p$, 
the holonomy of a periodic orbit along $\mathcal{C}_i$ is equal to $r_i V$, 
where $r_i$ is a rational number.
Hence $\sum_{j=0}^{q}\sum_{k=0}^{r(i,j)} v_{i,j,k} = r_i V$.

\begin{figure}[h!]
  \begin{center}
    \psfrag{riV}{$r_iV$}
    \psfrag{wi}{$w_i$}
    \psfrag{D1}{$\mathcal{D}_1$}
    \psfrag{D2}{$\mathcal{D}_2$}
    \psfrag{D3}{$\mathcal{D}_3$}
    \psfrag{D4}{$\mathcal{D}_4$}
    \psfrag{vi11}{$v_{i,1,1}$}
    \psfrag{vi21}{$v_{i,2,1}$}
    \psfrag{vi22}{$v_{i,2,2}$}
    \psfrag{vi23}{$v_{i,2,3}$}
    \psfrag{vi31}{$v_{i,3,1}$}
    \psfrag{vi41}{$v_{i,4,1}$}
    \includegraphics[width=.8 \linewidth]{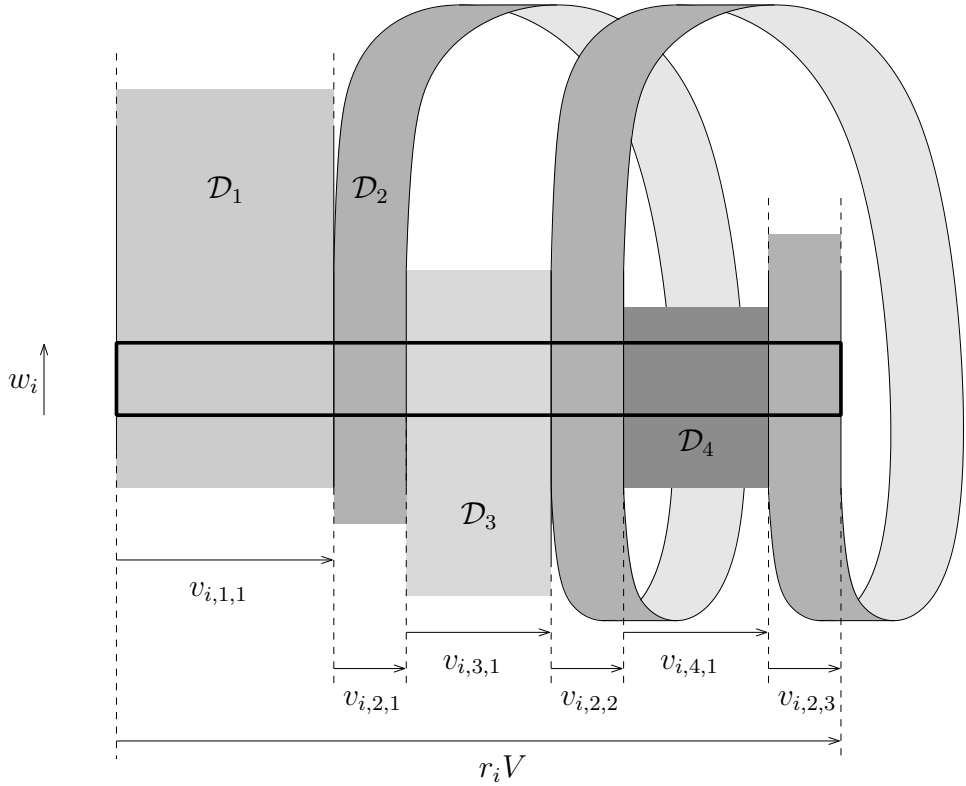}
    \caption{\label{pgm-cylindres.eps} Situation around $\mathcal{C}_i$.}
  \end{center}
\end{figure}

Since the closures of all the $P_{i,j,k}$'s form a cellullar decomposition of $\S$ 
(note that the singularities of $\S$ are located at the vertices of some $P_{i,j,k}$'s),
we have 
$$J(\S) = 
\sum_{i=0}^{p}\sum_{j=0}^{q}\sum_{k=0}^{r(i,j)} J(P_{i,j,k}) = 
\sum_{i=0}^{p} 2 r_i V \wedge w_i = 
2 V \wedge \sum_{i=0}^{p} r_i  w_i = 
V \wedge Z$$

where $Z$ is a vector colinear to $W$ that stands for $2 \sum_{i=0}^{p} r_i  w_i$.\\

The same computation holds if we replace $(V,W)$ by $(V',W)$, 
so we  have also  $J(\S) = V' \wedge Z'$ for some vector $Z'$ colinear to $W$.\\

So, we have $V \wedge Z =  V' \wedge Z'$ hence $(V,V',Z,Z')$ cannot be free over the rationals.
Hence, there exists $(a, b, c, d) \in \Q^4 \setminus \{(0,0,0,0)\}$ such that $aV+bV'+cZ+dZ'=0$.

Projecting on the $x$-axis along the $y$-axis, we have $aV+b\alpha V=0$ so $a=b=0$ as $\alpha\notin \Q$.
Therefore, we have $cZ+dZ'=0$ with $(c,d)\neq (0,0)$ and $(Z,Z')$ is not free over $\Q$.
So, there exists a rational number $q$ such that $Z=qZ'$ 
($Z$ and $Z'$ are different from $0$ since $J(\S)$ cannot be null).

Then $V \wedge qZ' =  V' \wedge Z'$, so $(qV-V') \wedge Z'  = (q-\alpha) V \wedge Z' = 0$.

Therefore, $(q-\alpha) V$ is colinear to $Z'$, a contradiction.

\item[Step 2: ] We prove that $hol(\P(\S))$ is a lattice of $\R^2$. 

Step 1 asserts that $dim_{\Q} (Vect_\Q(hol(P(\S))))=2$, 
and since translational holonomy is a morphism,
we have $Vect_\Q(hol(\P(\S)))\simeq \Q^2 \subset \R^2$.

Since $H_1(\S,\Z)$ is a free abelian group of finite type 
(it is isomorphic to $\Z^{2g}$ where $g$ denotes the genus of $\S$), 
so is $\P(\S)$ as a subgroup.

Any subgroup of finite type of $\Q^2$ is discrete in $\R^2$. 
Since there are periodic orbits in at least two directions, 
$hol(\P(\S))$ is a lattice of $\R^2$.
\end{description}\end{demo}

\begin{thm}\label{th-equivalence-homologie}
Let $\S$ be a translation surface whose homology is generated by the
periodic orbits of the geodesic flow on $\S$.
Then the following assertions are equivalent: 
\begin{enumerate}
  \item $\S$ is a torus branched covering
  \item $\S$ has the finite blocking property
  \item $\S$ is purely periodic
\end{enumerate}
\end{thm}

\begin{demo}
Since $(1)\Rightarrow(2)\Rightarrow(3)$ is already known
\cite{Monteil-pp}, it suffice to prove that 
a purely periodic translation surface $\S$ whose homology is generated by $P(\S)$ is a torus branched covering.
By Theorem \ref{th-hol=2}, $\Lambda \defi hol(H_1(\S,\Z)) = hol(\P(\S))$ is a lattice of $\R^2$.\\

Now, let us fix a point $x_0$ in $\S$. 
For each point $x$ in $\S$, we can chose a curve $\gamma_x$ from $x_0$ to $x$.
If $\gamma'_x$ is another such curve, then $hol(\gamma_x)-hol(\gamma'_x)\in \Lambda$.

Hence, the map $$\ba \S  \dans \R^2 / \Lambda \\ x \donne hol(\gamma_x) \mod \Lambda \ea$$ is well defined 
and realizes a branched covering from $\S$ to the torus $\R^2 /
\Lambda$ that preserves the translation structue (see
\cite{Zorich-branched}).
\end{demo}

This leads to the following general question:
\begin{center}{\em Is the homology of any translation surface generated by its periodic orbits?}\end{center}


Note that for our purpose, it is sufficient to prove that periodic orbits generate a subgroup of finite index of the holomogy to conclude.


\section{When is the homology generated by periodic orbits?}\label{section-generation-convex}
\subsection{For convex surfaces}

Let $\P$ be a simply connected polygon of $\R^2$ whose edges are grouped in pairs 
such that two edges in a same pair have the same length and direction but opposite orientation 
($\P$ gets the canonical orientation coming from $\R^2$ and transmits it to its boundary). 
Note that a ``geometrical'' side of $\P$ can be decomposed into more than one edge
(for exemple a rectangle can have more than four edges).
We say that $\P$ is a {\em pattern} of $\S$ if $\S$ can be obtained from $\P$ by identifying the edges in each pair by translation
(the singularities of $\S$ are therefore located at some vertices of $\P$).
The Veech construction of zippered rectangles ensures that any translation surface admits a pattern.
A {\em convex} translation surface is a translation surface that
admits a convex pattern \cite{Veech-hyperelliptic}.

\begin{prop}\label{th-convex}
Let $\S$ be a convex surface, then $P(\S)$ generates $H_1(\S,\Z)$.
\end{prop}

\begin{demo} 
Let $\P$ be a convex pattern of $\S$. Let us begin with the following observation:
\begin{description}
\item[Strolling round a singularity]
Let $\sigma$ be a singularity of $\S$ and let $v$ be a vertex of $\P$ that represents $\sigma$.
$v$ is surrounded by two consecutive edges in the counterclockwise order around $\P$: 
the first one is called the {\em left} edge and the second one is called the {\em right} edge around $v$ 
(see Figure \ref{bout-de-tour.eps}).

\begin{figure}[h!]
  \begin{center}
    \psfrag{v}{$v$}
    \psfrag{p}{$\P$}
    \psfrag{l}{left}    
    \psfrag{r}{right}
    \includegraphics[width=.5 \linewidth]{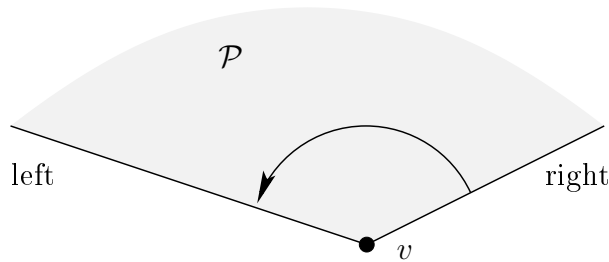}
    \caption{\label{bout-de-tour.eps} A piece of lap around $\sigma$.}
  \end{center}
\end{figure}

The left edge around $v$ is identified to the right edge around another vertex $v_1$ that also represents $\sigma$ in $\S$.
The left edge around $v_1$ is identified to the right edge around another vertex $v_2$ that also represents $\sigma$ in $\S$.
And so forth until we get a left edge around some $v_k$ that is identified
to the right edge around $v$ (see Figure \ref{left-right.eps}).

\begin{figure}[h!]
  \begin{center}
    \psfrag{vi}{$v_i$}
    \psfrag{vi1}{$v_{i+1}$}
    \psfrag{p}{$\P$}
    \includegraphics[width=.5 \linewidth]{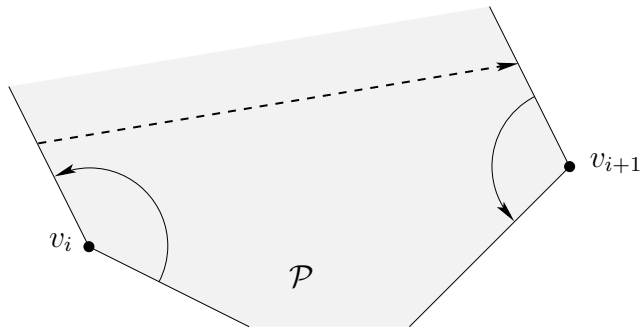}
    \caption{\label{left-right.eps} From the left edge around $v_i$ to the right edge around $v_{i+1}$.}
  \end{center}
\end{figure}

Wiewed from $\S$, we circled in the counterclockwise order around $\sigma$ (see Figure \ref{whole-tour.eps}).

\begin{figure}[h!]
  \begin{center}
    \psfrag{v}{$v_1$}
    \psfrag{p}{$\P$}
    \psfrag{l}{left}    
    \psfrag{r}{right}
    \includegraphics[width=.5 \linewidth]{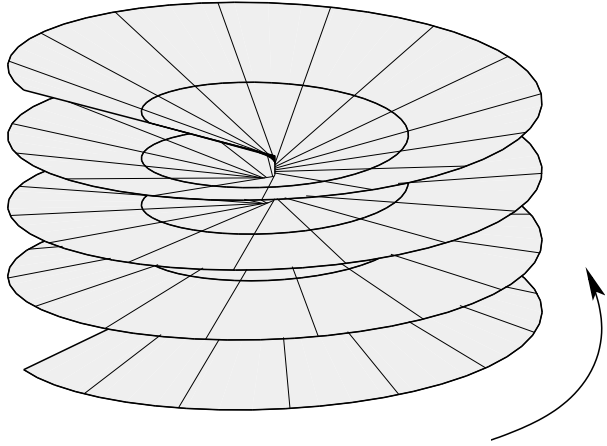}
    \caption{\label{whole-tour.eps} A whole lap around $\sigma$ in counterclockwise order.}
  \end{center}
\end{figure}

Hence, we passed around {\em every} vertex of $\P$ that corresponds to $\sigma$, 
and this stroll gives a cyclic order on the vertices of $\P$ that represent $\sigma$.

The observation is the following: the straight line from some $v_i$ to $v_{i+1}$ 
corresponds to an element of $\P(\S)$: 
indeed, viewed as a closed curve of $\S$, it is homologous to the 
closed curve in $\S$ associated to the straight line joining the middle of 
the left edge of $v_i$ to the right edge of $v_{i+1}$ that is a periodic orbit.

Hence, if $v$ and $v'$ are two vertices of $\P$ that represent the same singularity $\sigma$,
the straight line from $v$ to $v'$ in $\P$ corresponds in $\S$ to an element of $\P(\S)$, 
since it is homologous to the sum of the elementary 
closed curves associated to the straight lines from a $v_i$ to $v_{i+1}$ ($v'$ is one of the $v_j$'s).

\end{description}

Now, let $\gamma$ be a closed curve in $\S$.
Up to a tightening, we can assume that it is a closed geodesic.
If it is a periodic orbit, then we are done. 
Otherwise, it is a union of saddle connections $c_1$, $c_2$, \dots, $c_n$,
where $c_1$ goes from a singularity $\sigma_1$ to a singularity $\sigma_2$,  
$c_2$ goes from $\sigma_2$ to $\sigma_3$, and so forth, 
$c_n$ goes from $\sigma_n$ to $\sigma_1$ 
(the $c_i$'s belong to the relative homology and their sum belongs to the
absolute homlogy of $\S$).\\

Viewed in $\P$, $c_i$ is a finite union of segments $[a_{i,1},b_{i,1}]$, $[a_{i,2},b_{i,2}]$, \dots, $[a_{i,m_i},b_{i,m_i}]$, where
\begin{itemize}
  \item $a_{i,1}$ is a vertex of $\P$ that represents $\sigma_i$
  \item $b_{i,j}\sim a_{i,j+1}$ for $1\leq j \leq m_{i}-1$ (note that the endpoints of the segments belong to the boundary of $\P$) 
  \item $b_{i,m_i}$ is a vertex of $\P$ that represents $\sigma_{i+1}$ (with the convention that $n+1=1$)
\end{itemize}

For $1\leq i \leq n$ and $1 \leq j\leq m_{i}-1$, 
let $\lambda_{i,j}$ denotes the straight line in $\P$ between
$b_{i,j}$ and $a_{i,j+1}$: it corresponds in $\S$ to a periodic orbit
denoted by $\overline \lambda_{i,j}$.

For $1\leq i \leq n$, let $\mu_i$ denote the straight line in $\P$ between  $b_{i,m_i}$ and $a_{i+1,1}$ 
(with the convention that $n+1=1$): 
since $b_{i,m_i}$ and $a_{i+1,1}$ represent the same singularity $\sigma_{i+1}$, 
we observed in the beginning of the proof that $\mu_i$ corresponds in
$\S$ to an element $\overline \mu_i$ of $\P(\S)$.

Hence, viewed in $H_1(\S,\Z)$, the sum 
$$\gamma + \sum_{i=1}^n \sum_{j=1}^{m_{i}-1} \overline \lambda_{i,j} +
\sum_{i=1}^n \overline \mu_i$$
is homologous to zero: indeed, it is homologous to the image 
under the continuous projection $\P \rightarrow \S$ of a
closed path in the simply connected space $\P$ (see Figure \ref{chemin-convex.eps}).
Hence, $\gamma$ can be expressed as a sum of elements of $\P(\S)$ 
and the result is proved.

\begin{figure}[h!]
  \begin{center}
    \psfrag{a11}{$a_{1,1}$}
    \psfrag{a12}{$a_{1,2}$}
    \psfrag{b11}{$b_{1,1}$}
    \psfrag{b12}{$b_{1,2}$}
    \psfrag{a21}{$a_{2,1}$}
    \psfrag{a22}{$a_{2,2}$}
    \psfrag{a23}{$a_{2,3}$}
    \psfrag{a24}{$a_{2,4}$}
    \psfrag{b21}{$b_{2,1}$}
    \psfrag{b22}{$b_{2,2}$}
    \psfrag{b23}{$b_{2,3}$}
    \psfrag{b24}{$b_{2,4}$}
    \psfrag{m1}{$\mu_{1}$}
    \psfrag{m2}{$\mu_{2}$}
    \psfrag{l11}{$\lambda_{1,1}$}
    \psfrag{l21}{$\lambda_{2,1}$}
    \psfrag{l22}{$\lambda_{2,2}$}
    \psfrag{l23}{$\lambda_{2,3}$}
    \psfrag{r1}{represent $\sigma_1$}
    \psfrag{r2}{represent $\sigma_2$}
    \psfrag{rg}{represent $\gamma$}
    \includegraphics[width=.9 \linewidth]{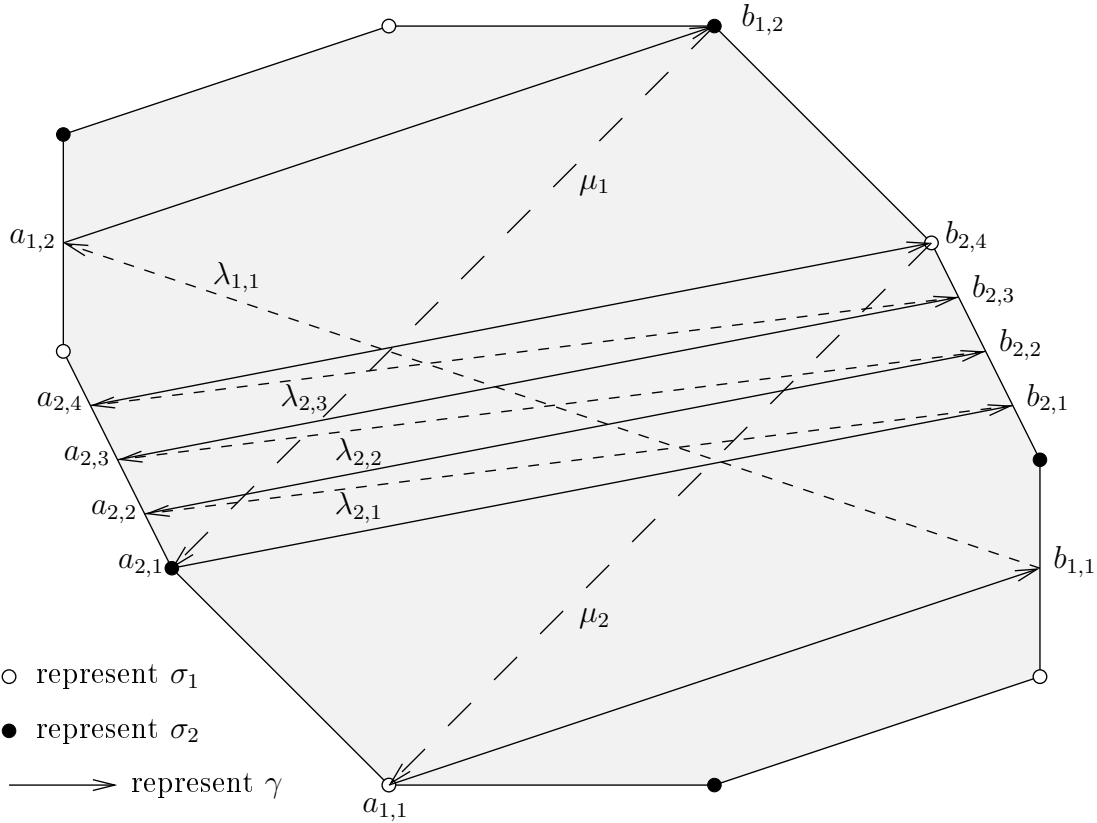}
    \caption{\label{chemin-convex.eps} We can add (dashed) elements of $\P(\S)$
        to $\gamma$ and obtain a closed curve homologous to zero in $\S$.}
  \end{center}
\end{figure}
 
\end{demo}

\begin{cor}
Any purely periodic convex translation surface is a torus branched covering.
\end{cor}

\begin{figure}[h!]
  \begin{center}
    \includegraphics[width=.5 \linewidth]{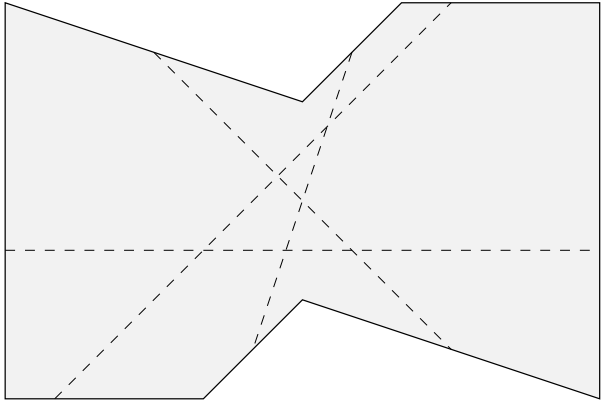}
    \caption{\label{face-to-face.eps} A face-to-face surface (identify parallel sides by translation).}
  \end{center}
\end{figure}

Note that we proved in fact the result for {\em face-to-face} translation surfaces i.e.
translation surfaces that 
admit a pattern $\P$ such that for each pair of identified edges $e\sim e'$,
there exists two identified points $x$ resp. $x'$ in $e$ resp. $e'$
such that the line between $x$ and $x'$ belongs to the interior of $\P$ (see Figure \ref{face-to-face.eps}).
Note that there exists non convex translation surfaces \cite{Veech-hyperelliptic},
but we do not know it there exists non face-to-face translation surfaces.

\subsection{On a dense open subset of full measure in any stratum}

\begin{prop}\label{th-generique}
In any stratum $\mathcal{H}_1(k_1,\dots,k_n)$, 
the set of translation surfaces whose homology is generated by the perodic orbits 
is an open dense subset of full measure. 
\end{prop}

\begin{demo}
Let $\C$ be a connected component of a normalized stratum and let $\mathcal{G}$ denotes the subset 
of surfaces whose homology is generated by the perodic orbits.

In \cite{KontsevichZorich}, Kontsevich and Zorich proved that there exists a translation surface 
that admits a cylinder decomposition in the horizontal direction with only one cylinder in $\C$.
Such surfaces are convex (the horizontal cylinder can be unrolled to a rectangle of $\R^2$) so $\mathcal{G}\neq \emptyset$.\\

Let $\S$ be an element of $\mathcal{G}$. 
So, there exists a finite number of periodic orbits $\gamma_1, \dots, \gamma_N$ in $\S$ that generate $H_1(\S,\Z)$.
Each $\gamma_i$ can be thickened to obtain a cylinder of positive height. 
Those cylinders survive under a small perturbation of $\S$,
so $\mathcal{G}$ is open in the stratum and therefore in $\C$.\\

Hence, $\mathcal{G}$ has positive measure, and since it is $SL(2,\R)$-invariant, it is of full measure in $\mathcal{G}$
(recall that the action of $SL(2,\R)$ is ergodic on $\C$). It is also dense in $\mathcal{G}$.\\

The result follows since the stratum is the finite union of its connected component.

\end{demo}

We have to notice that Zorich proved the same result using the asymptotic cycle, 
and looking to the evolution of the shape of Veech zippered rectangles of a generic surface along 
Rauzy inductions (see \cite{Zorich-private}).

We also obtained another result that is not exploited yet: if $\S$ is a
translation surface and $\mathcal{C}$ is a cylinder of $\S$ then there
exists $H>0$ such that for any $h\geq H$ the homology of the surface obtained by
replacing $\mathcal{C}$ by a cylinder of height $h$ and the same
perimeter as $\mathcal{C}$ is generated by its periodic orbits. 
Up to a renormalization of the area of those surfaces, the ``flow'' obtained 
when we let $h$ tend to infinity corresponds to a path from
$\S$ to the boundary of its normalized stratum in a manner that remains to study.


\section{Digression: Two directions do not suffice}

Recall that a surface is called {\em bouillabaisse} if its Veech group contains two transverse parabolic elements
(Veech surfaces are particular bouillabaisse surfaces).
Such surfaces admit cylinder decomposition in two directions (and the moduli are commensurable).
In \cite{Monteil-Veech}, we proved the equivalence of the three asssertions of theorem \ref{th-equivalence-homologie} for such surfaces. \\

This lets us think that it is sufficient to control cylinders in two directions to control the surface.
Indeed if $V$ and $W$ are two purely periodic directions (directions that admit a cylinder decomposition with commensurable perimeters), 
then the subgroup of $\R^2$ generated by the holonomy of periodic orbits of those two directions is a lattice and 
one might be tempted to apply the reasoning of theorem \ref{th-equivalence-homologie} to conclude.

We want to insist on the fact that two such directions are not generally sufficient to generate the homology, 
the following example should convince us:

\begin{figure}[h!]
  \begin{center}
    \psfrag{1}{$1$}
    \psfrag{2}{$2$}
    \psfrag{a}{$\alpha$}
    \includegraphics[width=.9 \linewidth]{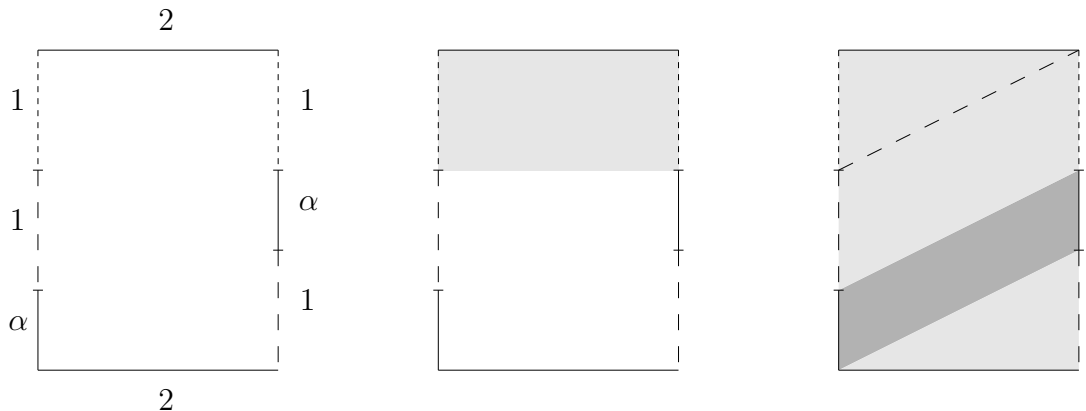}
    \caption{\label{2directions.eps} A strange surface; 
      parallel lines with the same style are identified by translation, 
      $\alpha\in \R\setminus \Q$.}
  \end{center}
\end{figure}

Indeed, this surface is purely periodic in the vertical direction and in the direction $(2,1)$ (third picture).
But there are in all three cylinders in those two directions, 
whereas the homology is a free abelian group on four generators. 
It is also not completely periodic, 
since the horizontal direction contains one cylinder and one minimal component (second picture).

Anyway, this surface is convex so its homology is generated by its periodic orbits, 
but one needs more than two directions to succeed.

\mathversion{bold}
\subsection{$J$-simple translation surfaces}
\mathversion{normal}

In the proof of theorem \ref{th-hol=2}, we used a particular property: 
let us call a translation surface {\em $J$-simple} 
if its $J$-invariant can be written $v\wedge w$ for some vectors $v$ and $w$ in $\R^2$.
Purely periodic translation surfaces are $J$-simple, 
but there exists non purely periodic $J$-simple translation surface, 
like the one described by Figure \ref{2directions.eps}.
Note that it is sufficient to have two transverse purely periodic directions to be $J$-simple. \\

Question: {\em is it necesssary to have two transverse purely periodic directions to be $J$-simple? 
can we describe $J$-simple translation surfaces?}

\subsection{Advertising}

We viewed in proposition \ref{th-generique} that is it not hard to obtain generic results about translation surfaces.
Most of those results do not hold everywhere, and to construct counter-examples is not an easy task.
For example, this section asked for $J$-simple translation surfaces
without two periodic directions, and surfaces whose homology cannot be
generated by periodic orbits in only two directions.
We want to propose here a tool to inventory all such constructions, 
an evolutive encyclopedia of {\em concrete} translation surfaces that
satisfy some particular properties. \\

Since the existing litterature and the folklore are quite fat already, 
and since the theory of translation surfaces is connected with many other fields
(ergodic theory, algebraic geometry, combinatorics, number theory, complex analysis,...), 
this objective can only be reached in an open and collborative way.\\

Technically, we set up a {\em wiki } i.e. a website where anyone can
freely create or modify any page using any web browser.
It is fully featured and lets the possibility to write {\LaTeX} formulas, to draw and add pictures, 
to add easily links to other pages,
to be automatically contacted when some new example is added, to follow the evolution of the web pages...
The website can be found at 
\begin{center} 
\texttt{http://ocarina.ath.cx/$\sim$titi/twiki/bin/view/WildSurfaces/WebHome} 
\end{center}

For example, we can inventory the
Veech's regular n-gons, 
a surface whose Veech group is $SL(2,\Z)$ but that is not a torus constructed by Herrlich and Schmithuesen,
McMullen's Veech $L$-shaped surfaces, 
some translation surfaces with non generic Siegel-Veech constants,
a $J$-simple surface without any purely periodic direction,
Veech obtuse triangular billiards of McBilliards, 
a non face-to-face translation surface,
a translation surface whose Hasudorff dimension of its non-uniquely ergodic directions is $1/2$ constructed by Cheung,
a translation surface whose Veech group is infinitely generated constructed by Hubert and Schmidt,
the Arnoux-Yoccoz surface whose Veech group contains an hyperbolic element but no parabolic element,
a surface that satisfies the Veech alternative but whose Veech group is not a lattice,
or whatever you want.
You can add new surfaces, but also 
ask for surfaces satisfying certain property, 
find new properties to an existing surface, 
simply add a comment, a new proof, a picture...\\

You will find on that website nothing more than what you will contribute.


\section*{Conclusion}

To sum up, the assertions
\begin{itemize}
  \item $\S$ is a torus branched covering (global geometric property)
  \item $\S$ has the finite blocking property (illuminatory property)
  \item $\S$ is purely periodic (dynamical property)
\end{itemize}
have been proved to be equivalent
\begin{itemize}
  \item for bouillabaisse surfaces, in particular Veech surfaces \cite{Monteil-Veech},
  \item in genus two \cite{Monteil-pp},
  \item for convex surfaces,
  \item on a dense open subset of full measure in any stratum.
\end{itemize}
Who is next?

\clearpage
\mathversion{bold}
\section*{Appendix: Why is the action of $SL(2,\R)$ on the moduli space of abelian differential so much studied?}

\vspace{2 cm}

\begin{figure}[h!]
  \begin{center}
   \psfrag{s}{}
   \includegraphics[width=.2 \linewidth]{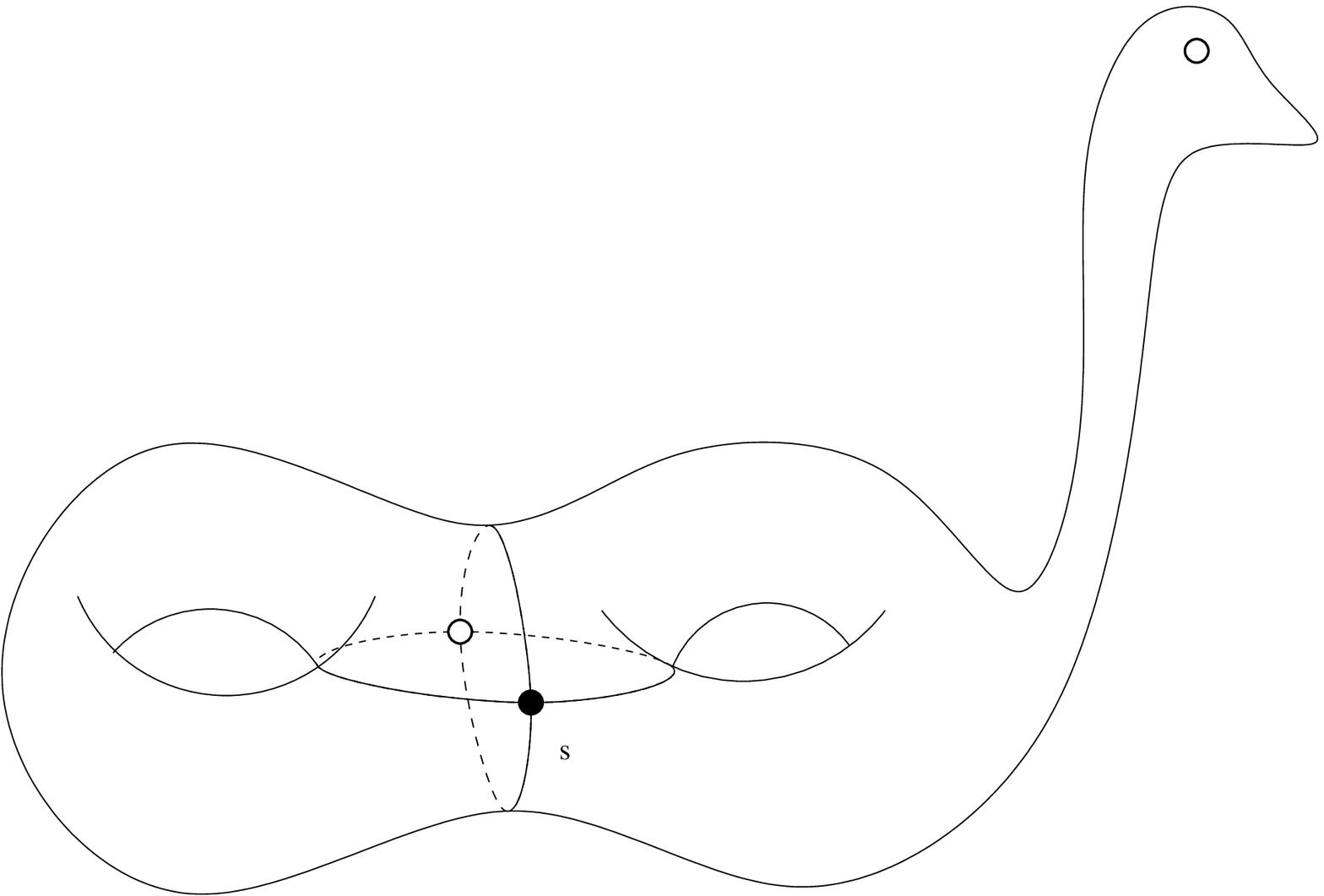}
  \end{center}
\end{figure}

\begin{figure}[h!]
  \begin{center}
  \psfrag{s}{$\sigma$}
  \includegraphics[width=.4 \linewidth]{coin.eps}
  \end{center}
\end{figure}

\begin{figure}[h!]
  \begin{center}
  \psfrag{s}{$\sigma$}
  \includegraphics[width=.6 \linewidth]{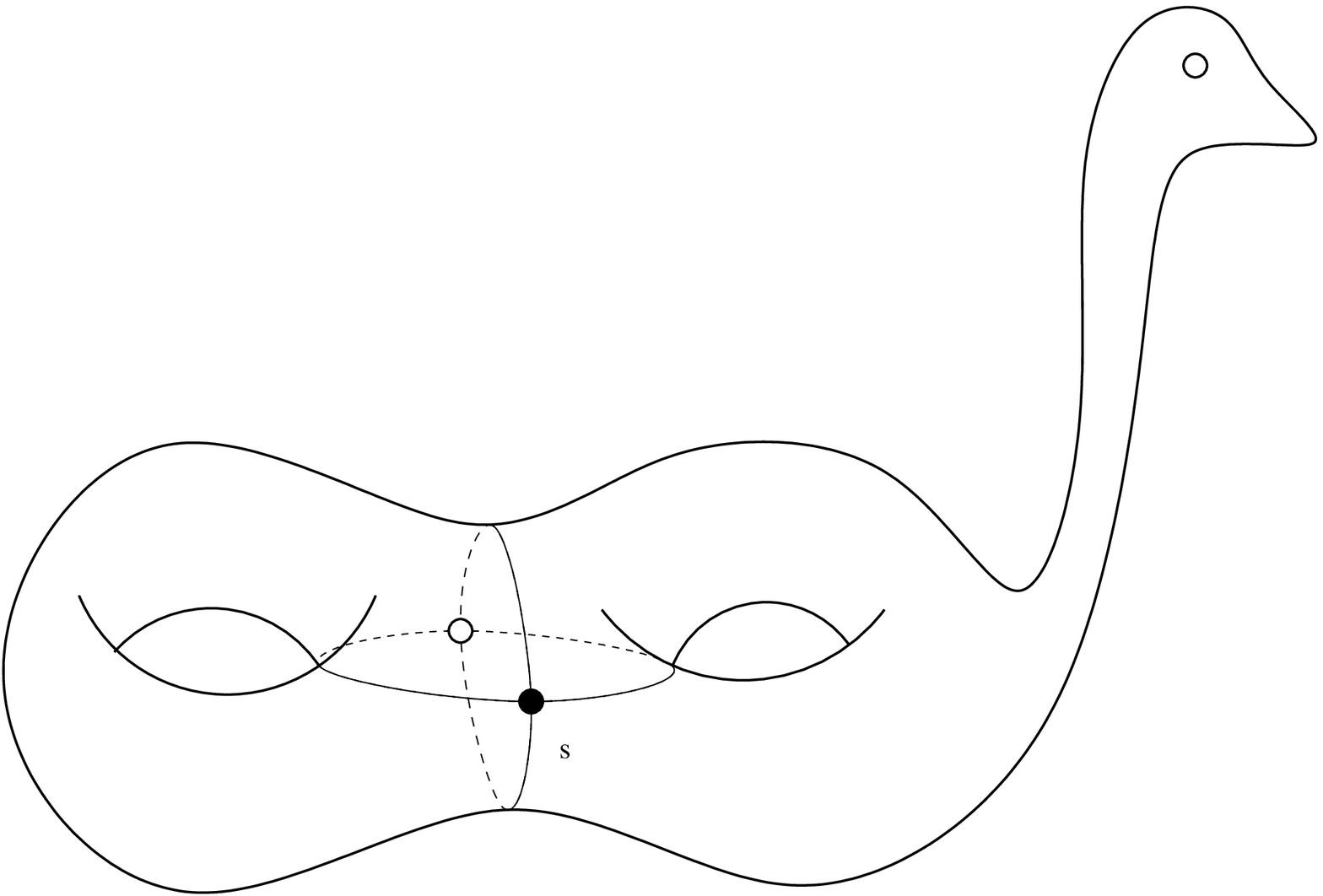}
  \end{center}
\end{figure}

\clearpage      

\begin{figure}[h!]
  \begin{center}
  \psfrag{s}{$\sigma$}
  \includegraphics[width=0.8 \linewidth]{coin2.eps}
  \end{center}
\end{figure}

\begin{figure}[h!]
  \begin{center}
  \psfrag{s}{\Huge \bf SPLASH!}
  \includegraphics[width=0.9 \linewidth]{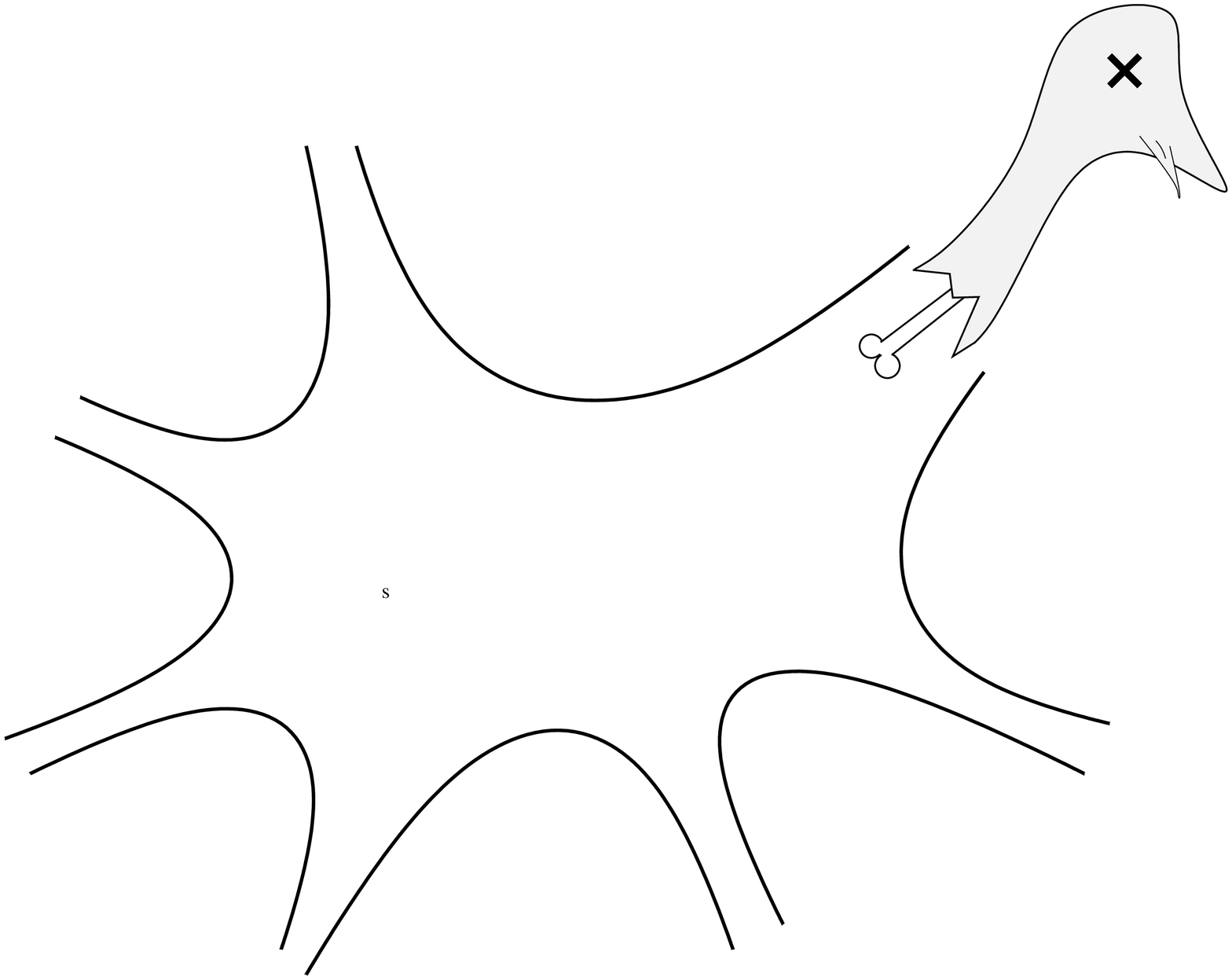}
  \end{center}
\end{figure}

\begin{center}
{\bf \Large Because the one of $GL(2,\R)$ 
is not optimal!}
\end{center}
\mathversion{normal}

\end{document}